\documentclass[12pt]{article}
\usepackage{amsmath,amssymb,amsthm}

\title{Orbifold Zhu theory associated to intertwining operators} 

\author{
  Hiroshi Yamauchi
  \vsb\\
  \it{\small Graduate School of Mathematics,}
  \\
  \it{\small University of Tsukuba, Ibaraki 305-8571, Japan}
  \\
  \it{\small e-mail:} \sf{\small  hirocci@math.tsukuba.ac.jp}
}

\date{}

\makeatletter
\@addtoreset{equation}{section}

\newcommand{\gs}[1]{\textbf{#1}}
\newcommand{\ds}{\displaystyle}
\newcommand{\fr}[2]{\frac{#1}{#2}}
\newcommand{\dfr}[2]{\dfrac{#1}{#2}}
\newcommand{\cd}{\cdot}
\newcommand{\cds}{\cdots}
\newcommand{\dsum}{\displaystyle \sum}

\newcommand{\simto}{\stackrel{\sim}{\to}}

\renewcommand{\l}{\left}
\renewcommand{\r}{\right}

\newcommand{\vsv}{\vspace{5mm}}
\newcommand{\vsb}{\vspace{2mm}}

\newcommand{\q}{\quad}

\newcommand{\maru}[1]{{\ooalign{\hfil#1\/\hfil\crcr
\raise.167ex\hbox{\mathhexbox20D}}}}

\newcommand{\ruby}[2]{%
 \leavevmode
 \setbox0=\hbox{#1}%
 \setbox1=\hbox{\tiny #2}%
 \ifdim\wd0>\wd1 \dimen0=\wd0 \else \dimen0=\wd1 \fi
 \hbox{%
   \kanjiskip=0pt plus 2fil
   \xkanjiskip=0pt plus 2fil
   \vbox{%
     \hbox to \dimen0{%
       \tiny \hfil#2\hfil}%
     \nointerlineskip
     \hbox to \dimen0{\mathstrut\hfil#1\hfil}}}}

\newcommand{\lfm}{\langle} 
\newcommand{\rfm}{\rangle} 
\newcommand{\la}{\langle}
\newcommand{\ra}{\rangle}

\newcommand{\abs}[1]{\lvert{#1}\rvert}

\newcommand{\tensor}{\otimes}

\newcommand{\Z}{\mathbb{Z}}
\newcommand{\C}{\mathbb{C}}

\newcommand{\N}{\mathbb{N}}
\newcommand{\Q}{\mathbb{Q}}

\newcommand{\res}{\mathrm{Res}}
\newcommand{\End}{\mathrm{End}}

\newcommand{\aut}{\mathrm{Aut}}
\newcommand{\wt}{\mathrm{wt}}
\newcommand{\tr}{\mathrm{tr}}
\renewcommand{\hom}{\mathrm{Hom}}

\newcommand{\SL}{\mathrm{SL}}

\newcommand{\pii}{\pi i}

\newcommand{\w}{\omega}
\newcommand{\vacuum}{\mathrm{1\hspace{-3.2pt}l}}
\newcommand{\vac}{\vacuum}

\newcommand{\im}{\mathrm{Im}}

\newcommand{\h}{\mathcal{H}}

\newcommand{\mat}{\begin{pmatrix} a & b \\ c & d \end{pmatrix}}

\newcommand{\pf}{\gs{Proof:}\q}

\theoremstyle{plain}
\newtheorem{thm}{Theorem}[section]
\newtheorem{prop}{Proposition}[section]
\newtheorem{lem}{Lemma}[section]

\theoremstyle{definition}
\newtheorem{df}{Definition}[section]

\theoremstyle{remark}
\newtheorem{rem}{Remark}[section]

\begin{document}

\baselineskip 6mm

\maketitle

\begin{abstract}
  We study intertwining operators among the twisted modules for
  rational VOAs and show a new modular invariance of the space 
  spanned by the trace functions associated to intertwining operators.
  Our result generalizes the results obtained by Zhu \cite{Z},
  Dong-Li-Mason \cite{DLM2} and Miyamoto \cite{M3}. 
\end{abstract}

\section{Introduction}

One of the main features of rational vertex operator algebras (VOAs
for short) is the modular invariance.
For a rational VOA $V$, Zhu \cite{Z} proved that the linear space
spanned by the trace functions $\tr_{W^i} o(a) q^{L(0)-c/24}$ 
defined on irreducible $V$-modules $W^i$ is invariant under the usual
action of the modular group $\SL_2 (\Z )$, where $o(a)$
denotes the grade-keeping operator (the zero-mode) of $a\in V$.
After his work, his result has been generalized in many directions.
Dong, Li and Mason generalized the theory above to the orbifold
theory in \cite{DLM2}. 
Miyamoto generalized the theory to involve intertwining operators 
in \cite{M3}.
In this paper, we give a new extension of the Zhu theory which
generalizes both the orbifold case \cite{DLM2} and the intertwining
operator case \cite{M3}.

Let $V$ be a VOA and $G$ a subgroup of $\aut (V)$.
A $V$-module $M$ is said to be {\it $G$-stable} if there exists a
group representation $\pi : G \to \End (M)$ (called the stabilizing
automorphism) such that $\pi (k) Y_M(a,z) v= 
Y_M( ka,z) \pi (k) v$ for all $a\in V$, $v\in M$ and $k\in G$.
Let $g$, $h$ be automorphisms on $V$ generating a finite abelian
subgroup $A$ of $\aut (V)$. 
Let $(U,\phi)$ be an $A$-stable untwisted $V$-module and $(W_g,\psi)$
an $h$-stable $g$-twisted $V$-module, where $\phi$ and $\psi$ are
stabilizing automorphisms on $U$ and $W_g$, respectively.
We introduce a notion of the $\phi (g)$-twisted intertwining
operator of type $U\times  W_g \to W_g$ (see Definition
\ref{twisted_intertwinnor}). 
Let $I(\cd,z)$ be such an intertwining operator.
Then one can consider the trace function associated to $I(\cd,z)$
defined as follows:
\begin{equation}\label{tr}
  S^I(u,\tau )=\tr_{W_g} z^{\wt u} I(u,z)\psi (h) q^{L(0)-c/24},
\end{equation}
where $u\in U$ and $q=e^{2\pii \tau}$.
In \cite{DLM2}, Dong, Li and Mason studied the trace functions
above in the case $U=V$, whereas Miyamoto considered them in the case
where $U$ is a general $V$-module and $A=1$ in \cite{M3}.
In each case, the space spanned by the trace functions has a modular
invariance property under certain conditions.
Therefore, it is a natural expectation that the linear space spanned
by the trace functions \eqref{tr} has a modular invariance
property for general $U$ and $A$.
The purpose of this paper is to prove that this is true.

This paper contains two main theorems.
The first main result is a theorem describing fusion rules among
untwisted and twisted modules in terms of the twisted Frenkel-Zhu's
bimodules (Theorem \ref{extended Li}). 
To prove the modular invariance in our case, we need a preliminary
result on the fusion rules.
So we will extend the theorem describing fusion rules among untwisted 
modules which is due to Frenkel-Zhu \cite{FZ} and Li
\cite{Li2} to the twisted case.
Using the extended Frenkel-Zhu-Li's theorem, we will prove that the
space spanned by the trace functions \eqref{tr} is invariant under the
action of the modular group, which is the second main theorem.

Let us explain our theorems more precisely.
Let $V$ be a VOA and $g$ an automorphism on $V$ of finite order.
Let $U$ be a $g$-stable $V$-module. We will construct an
$A_g(V)$-bimodule $A_g(U)$ as a quotient space of $U$ (see \cite{DLM1} 
for $A_g(V)$).
Using our bimodule $A_g(U)$, we can determine the fusion rules
(Theorem \ref{extended Li}).   
\vsb\\
\gs{Theorem 1.} 
{\it Assume that $V$ is $g$-rational.
Let $(U,\phi)$ be a $g$-stable $V$-module and $W^1_g$, $W^2_g$
irreducible $g$-twisted $V$-modules.
Then there is a linear isomorphism between
the linear space of the $\phi(g)$-twisted intertwining operators
of type $U\times W_g^1\to W_g^2$ and  
$\hom_{A_g(V)}\l( A_g(U)\tensor_{A_g(V)} W^1_g(0), W_g^2(0)\r)$, 
where $W_g^i(0)$ denotes the top levels of $W_g^i$ for $i=1,2$.}
\vsb\\
\indent 
Needless to say that, if we consider trivial automorphism, then our 
construction reduces to Frenkel-Zhu and Li's original one
(cf.{} \cite{FZ} \cite{Li2}).

Using the theorem above, we can show a new modular invariance. 
Let $g$, $h$ be automorphisms on $V$ generating a finite abelian
subgroup $A$ of $\aut (V)$.
Let $U$ be an $A$-stable $V$-module with a stabilizing automorphism
$\phi : A \to \End (U)$ and let $W_g$ be an $h$-stable $g$-twisted
$V$-module with a stabilizing automorphism $\psi (h)$.
Let $I(\cd,z)$ be a $\phi (g)$-twisted intertwining operator of type
$U\times W_g \to W_g$. We say that $I(\cd,z)$ is {\it $h$-stable} if
$\psi (h) I(u,z) \psi (h)^{-1}=I(\phi(h)u,z)$ holds for all $u\in U$.
Assume that $I(\cd,z)$ is $h$-stable. Then we consider the following
trace function:
\begin{equation}\label{tr2}
  S^I(u,\tau):=\tr_{W_g} z^{\wt u}I(u,z)\psi (h)^{-1} q^{L(0)-c/24}.
\end{equation}
Our main theorem is the following (Theorem \ref{conclusion}).
\vsb\\
\gs{Theorem 2.} 
{\it Assume that $V$ is $g$-rational. Let $(U,\phi)$ be an $A$-stable
$V$-module. Assume that $U$ is $C_{[2,0]}^A$-cofinite (see 
Definition \ref{C_{[2,0]}}).
Then the trace functions $S^I(u,\tau)$ defined by \eqref{tr2}
converge to holomorphic functions on the upper half plane.
Denote by $\mathcal{C}_1(U;\phi;(g,h))$ the linear space spanned by
the trace functions $S^I(\cd,\tau)$ associated to $h$-stable
$\phi(g)$-twisted intertwining operators $I(\cd,z)$ of type $U\times
W_g\to W_g$, where $W_g$ runs over $h$-stable $g$-twisted $V$-modules.
For $\gamma=\begin{pmatrix} a & b\\  f & d\end{pmatrix} \in \SL_2(\Z)$
with entries satisfying $a\equiv d\equiv 1 \mod \mathrm{lcm}
(\abs{g},\abs{h})$, $b\equiv 0\mod \abs{g}$ and $f\equiv 0 \mod
\abs{h}$, define its action on $S^I(\cd,\tau)$ as follows:
$$
  S^I|\gamma(u,\tau):= (f\tau +d)^{-k}S^I(u,\gamma \tau),
$$
where $u\in U_{[k]}=\{ u\in U\mid L[0]u=ku\}$.
Then $\gamma$ defines a linear isomorphism on
$\mathcal{C}_1(U;\phi;(g,h))$.}
\vsb\\
\indent
The argument used in the proof of Theorem 2 is almost the same as
those in \cite{Z}, \cite{DLM2} and \cite{M3}.
But it is worth pointing out the following remarkable difference.
In \cite{Z}, \cite{DLM2} and \cite{M3}, traces are taken over each
irreducible modules. 
But in our case we have to handle their conjugate under the
automorphism at the same time so that traces are taken over not only
irreducible modules.
This is the most difficult point of this paper.

The organization of this paper is as follows.
In Section 2 we recall basic definitions of the twisted theory of
VOAs that are needed in later sections.
The notion of the $\phi(g)$-twisted intertwining operators are also
given.
In Section 3, after introducing the twisted Frenkel-Zhu's bimodules, 
we extend Frenkel-Zhu and Li's theorem on fusion rules to the
twisted case.
In Section 4 we show the convergence of the trace functions defined as 
\eqref{tr}.
In Section 5 we prove the modular invariance of the space spanned by
the trace functions.

\section{Basic definition} 

The most notation are adopted from \cite{DLM2} and \cite{M3}.
Let $V$ be a VOA and $g$ an automorphism of $V$ of finite order
$\abs{g}$.
In this paper, we treat only $g$-rational VOAs.

\begin{df}\label{rational}
  A VOA $V$ is called $g$-rational if all admissible $g$-twisted
  $V$-modules are completely reducible.  
\end{df}

Let $(W_g,Y_g)$ be a $g$-twisted $V$-module. 
Take any $h\in \aut(V)$. 
Then one can find a natural $g^h=h^{-1}gh$-twisted $V$-module
structure on a vector space $W_g$. 
Define $g^h$-twisted vertex operator $(Y_g)^h$ on $W_g$ by
\begin{equation}\label{h-twist}
  (Y_g)^h(a,z):=Y_g(ha,z).
\end{equation}
Then $\l(W_g,(Y_g)^h\r)$ has a $g^h$-twisted $V$-module structure
(cf. \cite{DLM2}). 
We denote it by $W_g\circ h$.
In this paper, we consider the case where $g$ and $h$ generate a
finite abelian subgroup in $\aut (V)$.
A $g$-twisted $V$-module $W_g$ is called {\it $h$-stable} if 
$W_g\circ h$ is isomorphic to $W_g$. 
That is, there exists a linear isomorphism $\psi (h)$ on $W_g$
such that 
$$
  Y_g(ha,z)=\psi (h)Y_g(a,z)\psi (h)^{-1}.
$$ 
We call $\psi (h)$ {\it a stabilizing automorphism}.
Note that the stabilizing automorphism $\psi (h)$ is not
uniquely determined since non-zero scalar multiplications are allowed.

Let $U$ be a $g$-stable untwisted $V$-module admitting a weight
space decomposition $U=\oplus_{n\in\N} U_{n+h_0}$ and let $W_g^1$,
$W_g^2$ be irreducible $g$-twisted $V$-modules admitting weight space
decompositions $W_g^i=\oplus_{n\in\fr{1}{\abs{g}}\N} (W_g^i)_{n+h_i}$
for $i=1,2$, respectively, where we set $X_s:=\{ v\in X \mid L(0)
v=sv\}$ for a $V$-module $X$. 
One can find the general definition of intertwining operators among
twisted modules in \cite{DLM3} or \cite{X}, but we adopt a modified
definition as shown below.
Recall the definition of vertex operators on $g$-twisted $V$-modules.
Under the action of $g$, $V$ decomposes as follows:
$$
  V=V^0\oplus V^1\oplus \cds \oplus V^{\abs{g}-1},
$$
where $V^r:=\{ a\in V \mid ga=e^{2\pii r/\abs{g}}a\}$.
Let $M_g$ be an irreducible admissible $g$-twisted $V$-module.
By Theorem 8.1 (c) of \cite{DLM1}, one has a weight space
decomposition $M_g=\oplus_{n\in \fr{1}{\abs{g}}\N}(M_g)_{k+n}$, where
$k$ is a suitable rational number.
Setting $M_g^j:=\oplus_{n\in\N}(M_g)_{k+n-j/\abs{g}}$ for $0\leq j\leq 
\abs{g}-1$, we find that the vertex operators on $M_g$ provide the
following fusion rules among $V^0$-modules $V^i$ and $M_g^j$:
\begin{equation}\label{mane}
  V^i\times M^j_g=M^{i+j}_g,
\end{equation}
where by abuse of symbols we denote both an integer between $0$ and 
$\abs{g}-1$ and its residue class mod $\abs{g}$ by the same letter.
To introduce fusion rules similar to \eqref{mane},
we assume that $U$ has the following decomposition
\begin{equation}\label{decompo}
  U=U^0\oplus U^1\oplus \cds \oplus U^{\abs{g}-1}
\end{equation}
as a $V^0$-module and vertex operators on $U$ provide the following
fusion rules for $V^0$-modules:
\begin{equation}\label{fusion_assumption}
  V^i\times U^j= U^{i+j}.
\end{equation}
Based on the decomposition and fusion rules above, we define an
intertwining operator $I(\cd,z)$ of type $U\times W_g^1\to W_g^2$ as
follows.

\begin{df}\label{twisted_intertwinnor}
  Let $U=\oplus_{n\in\N} U_{n+h_0}$ be a $g$-stable $V$-module
  which has a decomposition \eqref{decompo} with fusion rules
  \eqref{fusion_assumption} and let  
  $W_g^i=\oplus_{n\in\fr{1}{\abs{g}}\N} (W_g^i)_{n+h_i}$, $i=1,2$, be
  irreducible $g$-twisted $V$-modules.
  An intertwining operator of type $U\times W_g^1\to W_g^2$ is a
  linear map
  $$
    I(\cd,z) : u\in U\mapsto
    I(u,z)=\sum_{n\in\C}u_nz^{-n-1}\in \hom \l(W_g^1, W_g^2\r)\{ z\} 
  $$
  with the following properties;
  \vsb\\
  $1^\circ$ \q for $u\in U^r$,
    $$
      I(u,z)\in \hom \l(W_g^1,W_g^2\r) 
      [[z,z^{-1}]]z^{-\fr{r}{\abs{g}}-h},
    $$
  where  $h=h_0+h_1-h_2$;
  \vsb\\
  $2^\circ$ \q for any $l\in\C$, there exists an $N\in \N$ such that
    $u_{l+n} w^1=0$ for all $n\geq N$;
  \vsb\\ 
  $3^\circ$ \q $L(-1)$-derivation : $I(L(-1)u,z)=\dfr{d}{dz}I(u,z)$;
  \vsb\\
  $4^\circ$ \q for $a\in V^r$ and $u\in U$, the following twisted
  Jacobi identity holds:  
    $$
     \begin{array}{c}
        z_0^{-1}\delta\l( \dfr{z_1-z_2}{z_0}\r)
          Y_{W_g^2}(a,z_1)I(u,z_2)
        -z_0^{-1}\delta\l( \dfr{-z_2+z_1}{z_0}\r) 
          I(u,z_2)Y_{W_g^1}(a,z_1)
        \vsb\\
        =z_2^{-1}\delta\l(\dfr{z_1-z_0}{z_2}\r) 
         \l(\dfr{z_1-z_0}{z_2}\r)^{-\fr{r}{\abs{g}}}
         I\l( Y_U(a,z_0)u,z_2\r) .
      \end{array}
    $$
\end{df}

Note that our assumptions \eqref{decompo} and
\eqref{fusion_assumption} are equivalent to the existence of a 
$g$-stabilizing automorphism $\phi (g)$ on $U$
such that the decomposition \eqref{decompo} coincides with the
eigenspace decomposition $U^r=\{ u\in U | \phi (g)u=e^{2\pii
r/\abs{g}}u\}$.  
Because our definition of intertwining operators is depend on the
choice of the decomposition of $U$, which owe to the choice of
$\phi (g)$, we call intertwining operators of type $U\times W_g^1\to
W_g^2$ defined as above ``{\it $\phi (g)$-twisted}'' to emphasize
the choice of the stabilizing automorphism $\phi (g)$.

The following two properties are direct consequences of the twisted
Jacobi identity;
\vsb\\ 
(1)\q Commutativity: for sufficiently large $N\in \N$, we have
$$
  (z_1-z_2)^N\{ Y_{W_g^2}(a,z_1)I(u,z_2)-I(u,z_2)Y_{W_g^1}(a,z_1)\}
  =0. 
$$
(2)\q Twisted Associativity: for $a\in V^r$ and $w^1\in W_g^1$, let
$k$ be a positive integer such that
$z^{k+\fr{r}{\abs{g}}}Y_{W_g^1}(a,z)w^1 \in W_g^1[[z]]$. Then
$$
  (z_0+z_2)^{k+\fr{r}{\abs{g}}} Y_{W_g^2}(a,z_0+z_2) I(u,z_2) w^1
  =(z_2+z_0)^{k+\fr{r}{\abs{g}}} I\l( Y_U(a,z_0)u,z_2\r) w^1.
$$
Conversely, we can prove the twisted Jacobi identity from these two
conditions (see \cite{Li1}).

\section{A generalization of Frenkel-Zhu-Li's theorem} 

In this section we construct an $A_g(V)$-bimodule $A_g(U)$ from a
$g$-stable $V$-module $U$ and show a theorem describing the fusion
rules of type $\binom{W_g^2}{U\ W_g^1}$ in terms of our bimodule
$A_g(U)$. 
These results are basically generalizations of the untwisted case.

\subsection{Associative algebra $A_g(V)$ and its bimodule $A_g(U)$}

Let $g$ be an automorphism of $V$ of finite order $\abs{g}$ and
$U$ a $g$-stable $V$-module with a stabilizing automorphism $\phi
(g)$. 
Recall the associative algebra $A_g(V)$ introduced in \cite{DLM1},
which is a quotient space of $V$ factored by a subspace 
$O_g(V)= \la a\circ_g b \mid a,b\in V\ra$. 
Let $a\in V^r=\{ v\in V \mid gv=e^{2\pii r/\abs{g}}v\}$.
We define a left action of $a$ on $U$ as follows:
\begin{equation}
  a \circ_gu 
  := \res_z \dfr{(1+z)^{\wt a+\fr{r}{\abs{g}}-1+\delta_{r,0}}}
     {z^{1+\delta_{r,0}}}Y(a,z)u.
\end{equation}
Using this product we define
$$
  O_g(U):=\lfm\, a\circ_g u\ |\ a\in V,\ u\in U\rfm 
$$
and set
$$
  A_g(U):=U^0/\l( O_g(U)\cap U^0\r) ,
$$
where $U^r:=\{ u\in U\, |\ \phi (g)u=e^{2\pii r/\abs{g}}u\}$.
We define a left action and a right action of $A_g(V)$ on this space.
By definition, $O_g(V)\supset \oplus_{i=1}^{\abs{g}-1} V^r$ so that
$A_g(V)$ is isomorphic to a quotient space of $V^0$. 
Hence we only need to consider the actions of $a\in V^0$ on $A_g(U)$.
For $a\in V^0$ we define left and right actions on $U$ as follows:
\begin{align}
  & a*_gu :=\res_z\dfr{(1+z)^{\wt a}}{z}Y(a,z)u,
  \vsb\\
  & u*_ga :=\res_z\dfr{(1+z)^{\wt a-1}}{z}Y(a,z)u .
\end{align}

When $g$ is trivial, the following theorem holds.

\begin{thm}\label{bimodule} 
  (\cite{FZ})\
  $A_1(U)$ is an $A_1(V)$-bimodule with respect to the actions $*_1$.
\end{thm}

As an extension of untwisted case we will show

\begin{thm}\label{twisted bimodule}
  $A_g(U)$ is an $A_g(V)$-bimodule with respect to the actions
  $*_g$.
\end{thm}

\pf 
We divide the proof into several lemmas.
The following lemma will be used frequently.
\begin{lem}\label{3.1}
  Let $a\in V^r$ be homogeneous and let $m\geq n\geq 0$. Then
  $$
    \res_z\dfr{(1+z)^{\wt
    a-1+\delta_{r,0}+\fr{r}{\abs{g}}+n}}{z^{m+1+\delta_{r,0}}} 
    Y(a,z)u \in O_g(U).
  $$
\end{lem}

\pf The proof is similar to that of Lemma 2.1.2. of \cite{Z}.
\qed
\vsb\\
We have to show  that $O_g(U)$ is closed under the actions of
$V^0$, the subspace $O_g(V)$ acts trivially on $U^0/O_g(U)\cap U^0$
and the associativity of the actions on $A_g(U)$.
By definition, we have $O_g(U)\cap U^0=O(U^0)+\sum_{r\ne 0}V^r\circ_g 
U^{\abs{g}-r}$, where $O(U^0)$ denotes the kernel of the (untwisted)
Frenkel-Zhu bimodule $A(U^0)$ associated to an untwisted $V^0$-module
$U^0$. 
By Theorem \ref{bimodule} we have $V^0*_g O(U^0)\subset
O(U^0)$ and $O(U^0)*_g V^0\subset O(U^0)$. 
Furthermore, for $a,b\in V^0$, we also have the following by Theorem
\ref{bimodule}.
$$
\begin{array}{c}
  a*_g O(U^0)\subset O(U^0),
  \vsb\\
  O(U^0)*_g a\subset O(U^0),
  \vsb\\
  a*_g( b*_g U^0) -(a*_gb)*_g U^0\subset O(U^0),
  \vsb\\
  ( U^0 *_g a) *_g b-U^0*_g(a*_gb)\subset O(U^0),
  \vsb\\
  ( a*_g U^0) *_gb-a*_g( U^0*_g b)\subset O(U^0).
\end{array}
$$
Now let $I=\sum_{r\ne 0}V^r\circ_g U^{\abs{g}-r}$. We should show
that $I$ is closed under the actions of $V^0$.

\begin{lem}
  Let $r\ne 0$, $a\in V^r$, $b\in V^0$ and $u\in U^{\abs{g}-r}$.
  Then both  $b*_g (a\circ_gu)$ and $(a\circ_g u)*_g b$ are contained
  in $O_g(U)\cap U^0$. 
\end{lem}

\pf The proof is similar to that of Proposition 2.3 of \cite{DLM1}.  
\qed

Next we will show an associative algebra $A_g(V)$ acts on  $A_g(U)$. 
In other words, we will show that
$$
\begin{array}{lll}
  \l( O_g(V)\cap V^0\r)*_g U^0 \subset O(U^0)+I
  &
  \text{and} 
  &
  U^0*_g \l( O_g(V^0)\cap V^0\r)\subset O(U^0)+I.
\end{array}
$$
By definition, we have $O_g(V)\cap V^0 = O(V^0) + \sum_{r\ne 0} V^r
\circ_g V^{\abs{g}-r}$ and by Theorem \ref{bimodule} we have
$O(V^0) *_g U^0 \subset O(U^0)$ and $U^0 *_g O(V^0) \subset O(U^0)$.
So it remains to show that $(V^r \circ V^{\abs{g}-r}) *_g U^0 \subset
O(U^0) +I$ and $U^0 *_g (V^r \circ_g V^{T-r}) \subset O(U^0) + I$ for
$r\ne 0$. 
Let $a\in V^r$, $r\ne0$, $b\in V^{\abs{g}-r}$ and $u\in U^0$. 
Then we have 
$$
\begin{array}{l}
  (a\circ_g b)*_g u=\res_w\dsum_{i=0}^\infty \binom{\wt
    a-1+\fr{r}{\abs{g}}}{i}\dfr{(1+w)^{\wt a+\wt
    b-i}}{w}Y(a_{i-1}b,w)u 
    \vsb\\
  =\res_w\res_z\dsum_{i=0}^\infty \binom{\wt
    a-1+\fr{r}{\abs{g}}}{i}\dfr{(1+w)^{\wt a+\wt b-i}}{w}
  \vsb\\
   \q \cd\l\{ (z-w)^{i-1}Y(a,z)Y(b,w)u-(-w+z)^{i-1}Y(b,w)Y(a,z)u\r\}
\end{array}
$$
$$
\begin{array}{l}
  =\res_w\res_z\dfr{(1+z)^{\wt a-1+\fr{r}{\abs{g}}}(1+w)^{\wt
    b+1-\fr{r}{\abs{g}}}}{w(z-w)}Y(a,z)Y(b,w)u
  \vsb\\
  \q -\res_w\res_z\dfr{(1+z)^{\wt a-1+\fr{r}{\abs{g}}}(1+w)^{\wt
    b+1-\fr{r}{\abs{g}}}}{w(-w+z)}Y(b,w)Y(a,z)u
  \vsb\\
  =\res_w\res_z\dsum_{i=0}^\infty \dfr{(1+z)^{\wt
    a-1+\fr{r}{\abs{g}}}}{z^{1+i}}\cd w^{i-1}(1+w)^{\wt b+1 -
    \fr{r}{\abs{g}}} Y(a,z)Y(b,w)u
  \vsb\\
  \q -\res_w\res_z\dsum_{i=0}^\infty (-1)z^i (1+z)^{\wt
    a-1+\fr{r}{\abs{g}}}\cd\dfr{(1+w)^{\wt b-1 +
    \fr{\abs{g}-r}{\abs{g}}+1}}{w^{2+i}} Y(b,w)Y(a,z)u
    \vsb\\
    \in O(U^0)+I.
\end{array}
$$
Similarly, we can show that $u *_g (a \circ_g b) \in O(U^0)+I$.
\qed

\subsection{Fusion rules}

Let $V$ be a rational VOA and let $W^i$, $i=1,2,3$, be irreducible
$V$-modules. 
The following formula is shown in \cite{FZ} and \cite{Li2}: 
$$
  \dim_\C \begin{pmatrix} W^3 \\ W^1\ W^2 \end{pmatrix} 
  =\dim_\C \hom_{A(V)}\l (A(W_1)\tensor_{A(V)}W_2(0), W_3(0)\r) .
$$
In this subsection we extend this formula to the twisted modules.
Let $U$ be a $g$-stable $V$-module with a stabilizing automorphism
$\phi (g)$ and let $W^1_g$, $W^2_g$ be irreducible $g$-twisted
$V$-modules. 
Assume that these modules admit weight space decompositions
$U=\oplus_{n\in \N}U_{n+h_0}$ and $W^i=\oplus_{n\in \fr{1}{\abs{g}}\N} 
(W^i_g)_{n+h_i}$, $i=1,2$, with some $h_j\in \Q$, $j=0,1,2$,
respectively. 
Set $U(n):=U_{n+h_0}$ and $W^i_g(n):=(W^i_g)_{n+h_i}$.
We may assume that {\it the top levels} $U(0)$, $W_g^i(0)$ of $U$,
$W_g^i$, respectively, are non-trivial.
We define $\deg u := \wt u -h_0$ and $\deg w^i:= \wt w^i-h_i$ for
$u\in U$ and $w^i\in W_g^i$, $i=1,2$, respectively.

By definition, a $\phi (g)$-twisted intertwining operator $I(\,\cd\
,z)$ of type $U\times W_g^1\to W_g^2$ has the form
$$
  I(u,z)=\sum_{n\in \Z +\fr{r}{\abs{g}}}u_{(n)}z^{-n-1-h}
$$
for $u\in U^r$, where $h=h_0+h_1-h_2$, and each coefficient 
$u_{(n)}$ maps a homogeneous subspace $W^1_g(m)$ of $W_g^1$ into a
homogeneous subspace $W^2_g(m+\deg u-n-1)$ of $W_g^2$. 
For homogeneous $u\in U^0$, we denote its weight-keeping operator (the
zero mode) by $o^I(u):=u_{(\deg u-1)}$. 
We also set $o^I(u)=0$ for $u\in U^r$, $r\ne 0$ and extend it linearly
on $U$. 
Similar to the untwisted case we have   

\begin{prop} 
  $o^I\in \hom_{A_g(V)}\l( A_g(U)\tensor_{A_g(V)}W^1_g(0),
  W^2_g(0)\r)$. 
\end{prop}

\pf Clearly $o^I$ is a linear map from $U^0\tensor_{\C}W^1_g(0)$
to $W^2_g(0)$. So we only need to show the following:
for $a\in V^0$ $u\in U^0$ and $w^1\in W_g^i(0)$,
\begin{eqnarray}
  & o^I\l( O_g(U)\r) w^1 =0;   \label{maru_1}\\
  & o^I\l( u*a\r) w^1=o^I\l( u\r)o(a)w^1; \label{maru_2}\\
  & \text{and}\q o^I\l( a*u\r) w^1=o(a)o^I\l( u\r) w^1. \label{maru3} 
\end{eqnarray}
First we consider \eqref{maru_1}.
By definition, $O_g(U)=\sum_{r=0}^{\abs{g}-1}V^r\circ_g
U^{\abs{g}-r}$, where $V^0\circ_g U^0=V^0\circ U^0$ (note that $\circ$
is the untwisted product) and we already know that $o^I\l
( V^0\circ_gU^0\r)=0$ in the untwisted theory (cf. \cite{FZ}). 
So it remains to show
$o^I\l( \sum_{r\ne 0} V^{r}\circ_{g}U^{\abs{g}-r}\r)=0$.
Since this proof is similar to that of Theorem 5.3 of \cite{DLM1}, we 
omit it.

Also \eqref{maru_2} and \eqref{maru3} are known for untwisted
intertwining operators and the proof is similar to the untwisted
case. For details, see \cite{FZ}.
\qed

Certain relations between $\hom_{A_g(V)}\l( A_g(U) \tensor_{A_g(V)} 
W^1_g(0), W^2_g(0)\r)$ and the space of intertwining operators are
known.

\begin{prop}(Proposition 2.10 of \cite{Li2})
  Assume that $W^i_g$ $(i=1,2)$ are  irreducible. Then the linear
  map  
  $$
  \begin{array}{llll}
    \pi : & I(\cd ,z) \in \dbinom{W_g^2}{U\ W_g^1} & \mapsto
    & o^I\in \hom_{A_g(V)}\l( A_g(U)\tensor_{A_g(V)}W^1_g(0), W^2_g(0)\r)
  \end{array}
  $$
  is injective.
\end{prop}

When $V$ is rational and $g=1$, the linear map $\pi$ defines an
isomorphism (Theorem 2.11 of \cite{Li2}). 
The following theorem is a generalization of Frenkel-Zhu-Li's
theorem. 

\begin{thm}\label{extended Li} 
  Let $V$ be a $g$-rational VOA and $(U,\phi (g))$ a $g$-stable 
  $V$-module. Let $W_g^i$, $i=1,2$, be irreducible $g$-twisted
  $V$-modules. Then 
  $$
    \dbinom{W_g^2}{U\ W_g^1}\simeq  
    \hom_{A_g(V)}\l( A_g(U)\tensor_{A_g(V)}W^1_g(0), W^2_g(0)\r) .
  $$
\end{thm}

\pf The argument is similar to that in the proof of Theorem 2.11 of 
\cite{Li2}. 
\qed

\begin{rem}\label{remark1}
  In the theorem above, we don't have to assume the irreducibility of
  $U$.
\end{rem}

\section{The space of trace functions} 

\subsection{The space of 1-point functions on the torus} 

In the following context, the following notation will be in force.
They are taken from \cite{DLM2}.
\vsb\\ 
(a)\ $V$ is a vertex operator algebra.
\vsb\\
(b)\ $g,h$ are automorphisms on $V$ generating a finite abelian
subgroup  $A=\la g,h\ra$ of $\aut (V)$. 
\vsb\\
(c)\ $g$ has order $\abs{g}$, $h$ has order $\abs{h}$ and $A$ has
  exponent $\mathrm{lcm}(\abs{g},\abs{h})$.
\vsb\\
(d)\ $\Gamma (g,h)$ is the subgroup of matrices $\mat$ in 
  $\SL_2(\Z )$ such that $b\equiv 0 \mod \abs{g}$, $c\equiv 0 \mod
  \abs{h}$ and $a\equiv d\equiv 1 \mod \mathrm{lcm}(\abs{g},\abs{h})$.
  \vsb\\
(e)\ $\mathcal{H}=\{ z\in \C | \im (z)>0\}$; the upper half plane. 
  \vsb\\
(f)\ $M(g,h)$ is the ring of holomorphic modular forms on $\Gamma
  (g,h)$.
  \vsb\\ 
(g)\ $U=\oplus_{n\in\N}U_{n+h}$ is an $A$-stable untwisted
  $V$-module\footnote{We need not assume the irreducibility of
  $U$. See also Remark \ref{remark1}.} 
  with a stabilizing automorphism $\phi : A\to
  \aut_\C (U)$ such that 
  $$
    Y_U (ka,z)=\phi (k) Y_U(a,z)\phi (k)^{-1}
  $$
  for all $k\in A$. 
  \vsb\\ 
(h)\ $U(g,h)=M(g,h)\tensor_{\C}U$.  
  \vsb\\
(i)\ $E_k(\tau )$, $P_k(\mu,\lambda,\tau )$, $Q_k(\mu,\lambda,\tau )$
  are defined as follows:
  $$
  \begin{array}{l}
    E_k(\tau )=\dfr{-B_k(0)}{k!} +\dfr{2}{(k-1)!}\dsum_{n=1}^\infty 
      \sigma_{k-1}(n)q^n ,\q k\geq 2,
    \vsb\\
    P_k(\mu, \lambda, z, q_\tau )=\dfr{1}{(k-1)!} \dsum_{n\in
      \fr{j}{M}+\Z}' \dfr{n^{k-1}q^n_z}{1-\lambda q_\tau^n} ,
      \q k\geq 1 ,
  \end{array}
  $$
  $$
  \begin{array}{l}
    Q_k(\mu,\lambda, q_\tau )=\dfr{1}{(k-1)!} \dsum_{n\geq 0}
      \dfr{\lambda (n+j/M)^{k-1}q_\tau^{n+j/M}}{1-\lambda
      q_\tau^{n+j/M}} 
    \vsb\\
    \hspace{3cm}
    + \dfr{(-1)^k}{(k-1)!}\dsum_{n\geq 1}\dfr{\lambda^{-1}
      (n-j/M)^{k-1} q_\tau^{n-j/M}}{1-\lambda^{-1} q_\tau^{n-j/M}}
    - \dfr{B_k(j/M)}{k!} ,\q k\geq 1 ,
  \end{array}
  $$
  and we also set $Q_0(\mu,\lambda,\tau )=-1$, where 
  $q_x$ denotes $e^{2\pii x}$, $\mu =e^{2\pii j/M}$, $\lambda =
  e^{2\pii l/N}$ for integers $j,l,M,N$ with $M,N>0$, the sign $\sum'$
  means omit the term $n=0$ if $(\mu,\lambda )=(1,1)$ and $B_k(x)$ are 
  the Bernoulli polynomials defined by the generating function
  $$
    \dfr{te^{tx}}{e^t-1}=\dsum_{k=0}^\infty B_k(x)\dfr{t^k}{k!}.
  $$
  For more informations, please refer to \cite{DLM2}.
  \vsb\\
We make it a rule to take 
$$
\begin{array}{lll}
  \mu^r=\abs{\mu}^re^{r(\arg \mu )i},
  & \text{where}\ \mu=\abs{\mu}e^{(\arg \mu )i},
  & -\pi < \arg \mu \leq\pi,
\end{array}
$$
for any $0\ne \mu\in\C$ and $r\in \Q$.
\vsb\\
For a VOA $(V,Y(\cd\, ,z),\vac, \w )$, one can define another VOA
structure on it.

\begin{df} (Theorem 4.2.1 of Zhu [Z]) 
  Let $(V,Y(\cd \, ,z),\vac,\w)$ be a VOA. For each homogeneous vector
  $a\in V$, the vertex operator 
  $$
    Y[a,z]:=Y(a,e^z-1)e^{z\wt (a)}=\sum_{n\in\Z}a_{[n]}z^{-n-1}\in
    \End (V)[[z,z^{-1}]] 
  $$
  gives a new VOA structure on $V$ with the same vacuum vector $\vac$
  and a new Virasoro vector $\tilde{\w}:=\w -\fr{c}{24}\vac$.
\end{df}
We denote the Virasoro operators of a new Virasoro vector $\tilde{\w}$ 
by $Y[\tilde{\w},z]$ $=$ $\sum_{n\in\Z}$ $L[n]$ $z^{-n-2}$.
Since $L[0]$ acts on $V$ semisimply, we can decompose $V$ into a
direct sum of eigenspaces for $L[0]$. Denote a new weight space
decomposition of $V$ by $V=\oplus_{n\in \N} V_{[n]}$, where
$V_{[n]}:=\{ v\in V \mid L[0]v=nv\}$. 
A $V$-module $U$ also admits a $L[0]$-weight space decomposition.
Denote by $U_{[k]}$ the $L[0]$-weight subspace of $U$ with weight $k$.
For a vector $u\in U_{[k]}$, we define its new weight by $\wt [u]
:=k$. 

Let $O(g,h)$ be an $M(g,h)$-submodule of $U(g,h)$ generated by the
following elements, where $a\in V$ satisfies $ga=\mu a$, $ha=\lambda
a$ and $u\in U$ satisfies $\phi (g) u =\mu' u$, $\phi (h) u=\lambda'u$
for some $\mu,\mu',\lambda,\lambda'\in\C$ ;
\begin{align}
  & a_{[0]}u,\ u\in U,\ (\mu,\lambda )=(1,1) , \label{O(g,h)1}
  \vsb\\
  & a_{[-2]}u+\sum_{k=2}^\infty (2k-1)E_{2k}(\tau )\tensor
  a_{[2k-2]}u,\ (\mu,\lambda )=(1,1) , \label{O(g,h)2}
  \vsb\\
  & u,\ (\mu',\lambda')\ne (1,1) , \label{O(g,h)3}
  \vsb\\
  & \sum_{k=0}^\infty Q_k(\mu,\lambda,\tau )\tensor a_{[k-1]}u,\ 
   (\mu,\lambda )\ne (1,1).
   \label{O(g,h)4}
\end{align}

We define the space of 1-point functions $\mathcal{C}_1(U;\phi;(g,h))$
to be the $\C$-linear space spanned by functions 
$$
  S : U(g,h)\times \h \to \C
$$
such that
\vsb\\
(C1)\ $S(u,\tau )$ is holomorphic in $\tau\in\h$ and for every $u\in
U(g,h)$;
  \vsb\\
(C2)\ $S(u,\tau )$ is $M(g,h)$-linear in the sense that $S$ is
$\C$-linear in $u$ and satisfies
\begin{equation}
  S(f\tensor u,\tau )=f(\tau )S(u,\tau )
\end{equation}
for $f\in M(g,h)$ and $u\in U$;
\vsb\\
(C3)\ $S(u,\tau )=0$ if $u\in O(g,h)$;
  \vsb\\
(C4)\ For $u\in U^A =\{ v\in U \mid \phi (k)v=v\ \text{for}\ k\in A\}$,
\begin{equation}\label{q5.8}
  S(L[-2]u,\tau )=\partial S(u,\tau )+\sum_{l=2}^\infty E_{2l}(\tau
  )S(L[2l-2]u,\tau ),
\end{equation}
where $\partial S$ is the operator which  is linear in $u$ and
satisfies
\begin{equation}\label{q5.9}
  \partial S(u,\tau )=\partial_kS(u,\tau )
  =\fr{1}{2\pii}\dfr{d}{d\tau}S(u,\tau )+kE_2(\tau )S(u,\tau )
\end{equation}
for $u\in U_{[k]}$.

Let $f(\tau )$ be a holomorphic function on $\h$.
We define the weight $k$ action of $\gamma =\mat\in\SL_2(\Z )$ by 
\begin{equation}
  (f|_k\gamma )(\tau ):=(c\tau +d)^{-k}f(\gamma \tau ) .
\end{equation}
We also set $(g,h)\gamma := (g^a h^c,g^b h^d)$.

The following theorem enables us to prove modular-invariance of trace
functions. 

\begin{thm}(\cite{DLM2} Theorem 5.4, Modular Invariance)
  \label{enable}
  For $S\in \mathcal{C}_1(U;\phi;(g,h))$ and $\gamma=\mat\in \SL_2(\Z
  )$, define
  \begin{equation}\label{q5.10}
    S|\gamma (u,\tau ):=S|_k\gamma (u,\tau )=(c\tau
    +d)^{-k}S(u,\gamma \tau )
  \end{equation}
  for $u\in U_{[k]}$ and extend linearly to $U$.
  Then $S|\gamma \in\mathcal{C}_1(U;\phi;(g,h)\gamma )$.
\end{thm}

In the rest of this subsection we assume that $U$ satisfies the
following finiteness condition which is weaker than the
$C_2$-cofiniteness.  

\begin{df}\label{C_{[2,0]}}
  Set $C_{[2,0]}^A (U):=\lfm a_{[-2]}u, b_{[0]}u \mid a\in V, 
  b\in V^A, u\in U\rfm$, where $V^A=\{\, x\in V \mid kx=x\ 
  \text{for}\ \text{all}\  k\in A\}$. 
  Then $U$ is said to be {\it $C_{[2,0]}^A$-cofinite} if
  $\dim (U/C_{[2,0]}^A(U))<\infty$. 
\end{df}

\begin{rem}
  The $C_{[2,0]}^1=C_{[2,0]}$-cofinite condition was
  introduced by Miyamoto in \cite{M3}. 
  It is clear from definitions that $C_2$-cofiniteness implies
  $C_{[2,0]}$-cofiniteness.
\end{rem}

Lemma 5.2 of \cite{DLM2} together with suitable modification leads the
following. 

\begin{lem}\label{qlem.5.2}
  Suppose $U$ is $C_{[2,0]}^A$-cofinite.
  Then $U(g,h)/O(g,h)$ is a finitely generated $M(g,h)$-module.
\end{lem}

\begin{lem}\label{qlem.5.3}
  Suppose $U$ is $C_{[2,0]}^A$-cofinite. Then for every
  $u\in U$ there exists $m\in\N$ and $r_i(\tau )\in M(g,h)$, $0\leq
  i\leq m-1$, such that
  \begin{equation}
    L[-2]^m u+\sum_{i=0}^{m-1}r_i(\tau )\tensor L[-2]^iu\in O(g,h) .
  \end{equation}
\end{lem}

\pf Since $U(g,h)/O(g,h)$ is a finitely generated $M(g,h)$-module and 
$M(g,h)$ is Noetherian, the submodule generated by $\{ L[-2]^iu \mid
i\geq 0\}$ is also a finitely generated submodule in
$U(g,h)/O(g,h)$. So the lemma follows. 
\qed
\vsb\\
We fix an element $S\in \mathcal{C}_1(U;\phi;(g,h))$.

\begin{lem}\label{qlem.6.1}
  Suppose $U$ is $C_{[2,0]}^A$-cofinite. Then for every $u\in U$ there
  exists some $m\in\N$ and $r_i(\tau )\in M(g,h)$ such that 
  \begin{equation}
    S(L[-2]^mu,\tau )+\sum_{i=0}^{m-1}r_i(\tau )S(L[-2]^iu,\tau )=0.
  \end{equation}
\end{lem}

\pf  Combine Lemma \ref{qlem.5.3} together with (C2) and (C3).
\qed
\vsb

Now $S\in \mathcal{C}_1(U;\phi;(g,h))$ has the same properties as that
in \cite{DLM2}. So we can apply exactly the same arguments as those in
\cite{DLM2} (Lemma 6.3-Lemma 6.8 in \cite{DLM2}) and hence we obtain
the following facts: 
for some $p\geq 0$, 
\begin{equation}\label{q6.12}
  S(u,\tau )=\sum_{i=0}^p (\log q_{\fr{1}{\abs{g}}})^i S_i(u,\tau ),
\end{equation}
where 
\begin{equation}\label{q6.13}
  S_i (u,\tau )=\sum_{j=1}^{b(i)}q^{\lambda_{ij}}S_{ij}(u,\tau ),
\end{equation}
\begin{equation}\label{q6.14}
  S_{ij}(u,\tau )=\sum_{n=0}^\infty a_{ijn}(u)q^{\fr{n}{\abs{g}}}
\end{equation}
are holomorphic on the upper half-plane, and
\begin{equation}\label{q6.15}
  \lambda_{ij_1}\ne \lambda_{ij_2}\q \l(\mathrm{mod}\
  \fr{1}{\abs{g}}\Z\r) 
\end{equation}
for $j_1\ne j_2$.

Hence we arrive at the following theorem which differs from
Theorem 6.5 in \cite{DLM2} only by changing of $C_2$-cofiniteness to
$C_{[2,0]}^A$-cofiniteness.  

\begin{thm}\label{qthm.6.5}
  Suppose $U$ is $C_{[2,0]}^A$-cofinite. 
  For every $u\in U$, the function
  $S(u,\tau )\in \mathcal{C}_1(U;\phi;(g,h))$ can be expressed in the
  form (\ref{q6.12})-(\ref{q6.15}). Moreover, $p$ is bounded
  independently of $u$.
\end{thm}

\subsection{$h$-conjugate intertwining operators} 

Let $A$, $U$ and $\phi$ be as previous.
By definition, the stabilizing automorphism $\phi$ gives a
representation of $\la g\ra$ on $U$.
Set $V^r:=\{ a\in V \mid ga=e^{2\pii r/\abs{g}}a\}$ and  
$U^r=\{ u\in U \mid \phi (g)u=e^{2\pii r/\abs{g}}u\}$ for 
$0\leq r\leq \abs{g}-1$.
Then $U=U^0\oplus U^1\oplus \cds \oplus U^{\abs{g}-1}$ and the vertex 
operators on $U$ give the following fusion rules for $V^0$-modules:
$$
  V^i\times U^j=U^{i+j}.
$$
Namely, the conditions \eqref{decompo} and \eqref{fusion_assumption}
are satisfied.
Therefore, a stabilizing automorphism $\phi : A=\la g,h\ra \to \aut_\C
(U)$ defines a $\la g\ra$-grading decomposition on $U$ and hence 
we can think of $\phi (g)$-twisted intertwining operators.

Let $W_g^i$, $i=1,2$, be  irreducible $g$-twisted $V$-modules.
Recall the $h$-conjugates $W_g^i\circ h$ of $W_g^i$ defined by
\eqref{h-twist}. By definition, there exist linear maps
$\psi_i (h):W_g^i\circ h \to W_g^i$ $(i=1,2)$ such that
\begin{equation}\label{h-trans}
  Y_{W_g^i}(ha,z)=\psi_i (h) Y_{W_g^i \circ h}(a,z)\psi_i (h)^{-1}.
\end{equation}
Assume a $\phi (g)$-twisted intertwining operator $I(\cd\, ,z)$
of type $U\times W_g^1\to W_g^2$ is given. 
We can define its $h$-conjugate intertwining operator 
$I^h(\cd\, ,z)$ of type $U\times W_g^1\circ h\to W_g^2\circ h$ as
follows:
\begin{equation}\label{h-conj}
  I^h(u,z):= \psi_2(h)^{-1} I\l( \phi (h)u,z\r)
  \psi_1(h). 
\end{equation}
Therefore, there exists a linear isomorphism 
$$
  \binom{W_g^2}{U\ W_g^1}
  \simeq \binom{W_g^2\circ h}{\ U \q W_g^1\circ h}
$$ 
given by an $h$-conjugation $I(\cd\, ,z) \mapsto I^h(\cd\, ,z)$.

\subsection{Fundamental $h$-stable sectors} 

Let $\mathcal{M}=\{ (W_g^\alpha ,Y^{\alpha})\}_{\alpha\in \Lambda}$ be 
the set of all inequivalent irreducible $g$-twisted $V$-modules. 
Since we have assumed that $V$ is $g$-rational, $\Lambda$ is a finite
set and every irreducible $g$-twisted $V$-module has a weight space 
decomposition such that all its homogeneous subspaces are of finite
dimension (cf. \cite{DLM1}).  
We consider a certain equivalent relation on $\mathcal{M}$.
Let $(W_g^\alpha,Y^\alpha )$ be an irreducible  $g$-twisted
$V$-module in $\mathcal{M}$. 
Then an $h$-conjugate $\l( W_g^\alpha ,(Y^\alpha )^h\r)$ is also an
irreducible $g$-twisted $V$-module and hence there exists a unique
$\beta \in \Lambda$ such that $\l( W_g^\alpha, (Y^\alpha )^h\r)$ is
isomorphic to $(W_g^\beta, Y^\beta )$. 
Therefore, by defining $h(\alpha )=\beta$ for each $\alpha\in \Lambda$
in this way, $h$ defines a permutation on the index set $\Lambda$ and
hence on $\mathcal{M}$.
We set $(W^\alpha,Y^\alpha)\circ h :=(W^{h(\alpha)}, Y^{h(\alpha)})$.
By definition, there exists unique linear isomorphisms $\psi_\alpha$
up to scalar multiplications such that
the following commutative diagrams hold for all $\alpha \in
\Lambda$: 
$$
\begin{array}{lccr}
   W_g^{h(\alpha )} & \xrightarrow{\ Y^{h(\alpha)} (a,z)\ } 
   & W_g^{h(\alpha )} 
  \vsb\\
  \downarrow\psi_\alpha & \circlearrowleft & \downarrow \psi_\alpha 
  \vsb\\
  W_g^\alpha & \xrightarrow{\ Y^\alpha (ha,z)\ } 
  & W_g^\alpha .
\end{array}
$$
We define an equivalent relation on $\mathcal{M}$ as follows.
Define $(W^\alpha_g, Y^\alpha )\sim (W^\beta_g, Y^\beta)$ if  
both $(W^\alpha_g, Y^\alpha )$ and $(W^\beta_g, Y^\beta )$ belong to
the same orbit under the action of $h$, or equivalently there exists
some $i\geq 0$ such that $(W^\alpha_g, Y^\alpha )\circ h^i\simeq 
(W^\beta_g,Y^\beta)$. 
Take an element $(W_g^\alpha,Y^\alpha)$ in $\mathcal{M}$ and fix it.
Assume its orbit $\{ (W_g^\alpha, Y^\alpha )\circ h^i\}$ ($i=0, 1, \dots, 
n-1$) is of length $n$. 
For simplicity we denote $(W_g^\alpha, Y^\alpha )\circ h^i$ by
$(W_g^i,Y^i)$. 
Then we have linear isomorphisms $\psi_i(h) : W_g \circ h^{i+1} \to
W_g \circ h^i$, $i=0,\cdots,n-1\, (\in \Z /n\Z)$, such that
\begin{equation}\label{Y^i}
  Y^i(ha,z)=\psi_i(h)Y^{i+1}(a,z)\psi_i(h)^{-1}.
\end{equation}
Set $\ds \bar{W}_g:=\bigoplus_{i=0}^{n-1} W_g^i$ and $\bar{\psi}(h):=
\psi_0(h)\oplus \cds \oplus \psi_{n-1}(h)$.
Then $\bar{W_g}$ with an automorphism $\bar{\psi}(h)$ is the minimal
$h$-stable $g$-twisted $V$-module in the sense that it contains $W_g^i$,
$i=0,1,\dots, n-1$, as submodules with multiplicity one. 
In this paper, an $h$-stable $g$-twisted $V$-module $\bar{W}_g$
constructed as above is of special importance so that we call the
modules constructed as above {\it fundamental $h$-stable $g$-twisted
$V$-modules}.

\subsection{Trace functions} 

Let $M_g$ be a fundamental $h$-stable $g$-twisted $V$-module with
an $h$-stabilizing automorphism $\psi (h)$.
Because $V$ is $g$-rational, there exists an irreducible $g$-twisted
$V$-submodule $W_g$ such that $M_g=\oplus_{i=0}^{n-1} W_g\circ h^i$,
where $n$ is the minimum positive integer such that $W_g\circ h^n$ 
is isomorphic to $W_g$.  
Let $I(\cd\, ,z)$ be a $\phi (g)$-twisted intertwining operator of
type $U\times M_g\to M_g$. 
Since $M_g$ is $h$-stable, an $h$-conjugate intertwining
operator $I^h(\cd\, ,z)=\psi (h)^{-1} I(\phi (h)\, \cd\, ,z)\psi (h)$
is also of type $U\times M_g\to M_g$.
Hence, $h$ acts on the space $\binom{M_g}{U\ M_g}$.
Since we need a nice relation between $\phi (h)$ and $\psi (h)$, we 
adopt the following definition.

\begin{df}
  A $\phi (g)$-twisted intertwining operator $I(\cd\, ,z)$ of type
  $\binom{M_g}{U\ M_g}$ is called {\it $h$-stable} if $\psi_h I(u
  ,z)\psi_h^{-1}=I(\phi (h)u ,z)$ holds for any $u\in U$.
\end{df}

The condition that $I(\cd,z)$ to be $h$-stable is equivalent to that
$I(\cd,z)$ is stable under the action of $h$ on $\binom{M_g}{U\
M_g}$. This is why we call ``stable''. 
Assume that $I(\cd,z)$ is $h$-stable. The main subject we study in the
rest of this paper is the trace functions defined as follows:
\begin{equation}\label{trace_function}
  S^I (u,q) := z^{\wt u} \tr_{|_{M_g}} I(u,z) \psi(h)^{-1}
  q^{L(0)-c/24}.
\end{equation}
We claim that the trace function \eqref{trace_function} associated to 
an $h$-stable $\phi (g)$-twisted intertwining operator $I(\cd,z)$
belongs to $\mathcal{C}_1(U;\phi;(g,h))$. 
We may present a rigorous proof here, but it will be a repetition of  
some complicated formal calculations which are the same as those in
\cite{Z}, \cite{DLM2} and \cite{M3} so that we omit details.
For details, please refer to Section 8 of \cite{DLM2} and Section 3 of 
\cite{M3}. Here is the consequence.

\begin{thm}\label{prev}
  Let $U$ be an $A$-stable untwisted $V$-module with an
  $A$-stabilizing automorphism $\phi$,
  and assume that $U$ is $C_{[2,0]}^A$-cofinite.
  Let $M_g$ be a fundamental $h$-stable $g$-twisted $V$-module with
  an $h$-stabilizing automorphism $\psi (h)$ and let $I(\cd\, ,z)$ be
  an $h$-stable $\phi (g)$-twisted intertwining operator of type
  $U\times M_g\to M_g$. 
  Then the trace function
  $$
    S^I (u,q) := z^{\wt u} \tr_{|_{M_g}} I(u,z) \psi(h)^{-1}
    q^{L(0)-c/24}, \q q=e^{2\pii \tau}
  $$
  gives an element of $\mathcal{C}_1(U;\phi;(g,h))$.
\end{thm}

By definition, a fundamental $h$-stable $g$-twisted $V$-module $M_g$
is a direct sum of irreducible modules.
So an $h$-stable $\phi (g)$-twisted intertwining operator $I(\cd,z)$
of type $U\times M_g\to M_g$ is also a direct sum of intertwining
operators among irreducible components.
Let $W_g$ be an irreducible submodule of $M_g$.
It is clear that on $W_g$, only an intertwining operator of type 
$U\times W_g\circ h \to W_g$ contributes to the trace function
\eqref{trace_function}.  
So we may assume that $I(\cd,z)$ is a direct sum of intertwining
operators of type $U\times W_g\circ h^{i+1} \to W_g\circ h^i$ with
$i\in \N$.
We call such an $h$-stable $\phi (g)$-twisted intertwining operator
{\it fundamental}.

\section{Modular invariance} 

In the previous section we showed that the trace function 
$S^I(\cd,\tau)$ associated to an $h$-stable intertwining operator
$I(\cd,z)$ belongs to $\mathcal{C}_1(U;\phi;(g,h))$. 
We will prove here that we can choose a basis of the linear space
$\mathcal{C}_1(U;\phi;(g,h))$ consisting of trace functions
$S^I(\cd,\tau)$. 
Namely, we shall show that the linear space spanned by $S^I(\cd,\tau)$
coincides with $\mathcal{C}_1(U;\phi;(g,h))$, and at the same time we
prove a new kind of modular invariance because the space
$\mathcal{C}_1(U;\phi;(g,h))$ is invariant under the action of the 
subgroup $\Gamma (g,h)$ of the modular group $\SL_2(\Z )$ by Theorem
\ref{enable}.

First, we need the following lemma.  

\begin{lem}\label{qlem.9.2} (\cite{DLM2} Lemma 9.2)\ 
  Suppose $a\in V$ satisfies $ga=\mu a$, $ha=\lambda a$ with 
  $\mu,\lambda\in\C$. Then the following hold:
  \vsb\\
  (i)\ Seen as an element of $U[[q_{\tau }]]$,
    the constant term of $a_{[-1]} u - \sum_{k=1}^\infty E_{2k}
    (\tau )$ $a_{[2k-1]} u$ is equal to $a*_gu-\fr{1}{2}a_{[0]}u$
    and the constant term of $a_{[-2]} u +\sum_{k=2}^\infty (2k-1)
    E_{2k}(\tau ) a_{[2k-2]}u$ is equal to $\fr{1}{12} a_{[0]}
    u+a\circ_gu$ if $\mu =\lambda =1$.
  \vsb\\
  (ii)\ The constant term of  $\sum_{k=0}^\infty
    Q_k(\mu,\lambda,q) a_{[k-1]}u$ is equal to $-a\circ_g u$ 
    if $\mu\ne 1$, and the constant term of $\sum_{k=0}^\infty
    Q_k(\mu,\lambda,q)a_{[k-1]}u$ is equal to $-a*_g u+
    \fr{1}{1-\lambda} a_{[0]}u$ if $\mu =1,\lambda\ne 1$. 
\end{lem}

Because $(L[-1]a)*_g u=a\circ_g u$ and $(L[-1]a)_{[0]}=0$, by 
replacing $a$ by $L[-1]a$ in the case where $\mu =1$ and $\lambda\ne
1$, we note that $-a\circ_g u$ also appears in the constant term in
that case. 
Thus, the lemma above insists that almost all generators of $O_g(U)$ 
appear in the constant terms of elements in $O(g,h)$. 

Let us study the relation between $\mathcal{C}_1(U;\phi;(g,h))$ and
our trace functions. Take an arbitrary function $S\in
\mathcal{C}_1(g,h)$. 
By Theorem \ref{qthm.6.5}, we know that $S$ is expressed as
$$
  S(u,\tau )=\sum_{i=0}^p(\log q_{\fr{1}{\abs{g}}})^iS_i(u,\tau )
$$
for a fixed $p$ and all $u\in U$ with each $S_i$ satisfying 
(\ref{q6.13})-(\ref{q6.15}).

We define a new operation  $\tilde{\w} *_{\tau}$ on $U$ by
$$
  \omega *_{\tau}u := L[-2]u
  - \dsum_{k=1}^\infty E_{2k}(\tau )L[2k-2]u 
$$
for $u\in U$. For convenience, we also write 
$$
  \tilde{\w}*_{\tau}S(u,\tau )=S(\tilde{\w}*_{\tau}u,\tau ) .
$$
In the following, $q_x$ means $e^{2\pii x}$.

\begin{lem} 
  We have
  \begin{align}
    & S_{ij}(a_{[0]}u,\tau )=0\q \text{for}\ \text{any}\ a\in V^A,
       u\in U, i,j,
       \label{m4.15}
    \vsb\\
    & S_{ij}(u_q,\tau )=0 \q 
      \text{for}\ \text{any}\ u_q\in O(g,h), i,j,
      \label{m4.16}
    \vsb\\
    & S_{pj}(\tilde{\w}*_{\tau}u,\tau )
      = \dfr{1}{\abs{g}}q_{\fr{1}{\abs{g}}} 
        \dfr{d}{dq_{\fr{1}{\abs{g}}}} 
        S_{pj}(u,\tau ) \q \text{for}\ u\in U,
        \label{m4.17}
  \end{align}
  \begin{align}
    & S_{i}(\tilde{\w}*_{\tau}u,\tau )
      = \dfr{1}{\abs{g}}\l\{ (i+1)S_{i+1}(u,\tau )
      + q_{\fr{1}{\abs{g}}} \dfr{d}{dq_{\fr{1}{\abs{g}}}} 
        S_i(u,\tau )\r\}
      \q \text{for}\ u\in U,
        \label{m4.18}
    \vsb\\
    &
    \text{and}\q
    \l(\tilde{\w}*_{\tau}-\dfr{1}{\abs{g}}q_{\fr{1}{\abs{g}}}
      \dfr{d}{dq_{\fr{1}{\abs{g}}}}
    \r)^N 
    \cd S_{ij}(u,\tau )=0 \q \text{for}\ N\gg 0.
       \label{m4.19}
  \end{align}
\end{lem}

\pf  Since $O(g,h)\in U[[q^{\fr{1}{\abs{g}}}]]$, (\ref{m4.15}) and 
(\ref{m4.16}) are trivial. 
Note 
$$
  \fr{1}{2\pii}\fr{d}{d\tau} = \dfr{1}{\abs{g}}
  q_{\fr{1}{\abs{g}}} \dfr{1}{dq_{\fr{1}{\abs{g}}}}
$$
so that from (C4) we have 
$$
  S(\tilde{\w}*_{\tau}u,\tau )
  = \dfr{1}{\abs{g}}\dsum_{i=0}^p\l\{ 
      i\l(\log q_{\fr{1}{\abs{g}}}\r)^{i-1}S_i(u,\tau )
    + \l(\log q_{\fr{1}{\abs{g}}}\r)^i
      q_{\fr{1}{\abs{g}}} \dfr{d}{dq_{\fr{1}{\abs{g}}}} S_i(u,\tau )
    \r\} .
$$
Therefore, we get (\ref{m4.17}) and (\ref{m4.18}). 
In particular, we find that
$$
  \l( \tilde{\w}*_{\tau} - \fr{1}{\abs{g}} q_{\fr{1}{\abs{g}}} 
  \fr{1}{q_{\fr{1}{\abs{g}}}}\r) S_i(u,\tau ) 
  = \fr{1}{\abs{g}} (i+1) S_{i+1}(u,\tau )
$$
and hence (\ref{m4.19}) also holds.
\qed
\vsv\\
We note that both $S_i(u,\tau )$ and the trace functions
$S^I(u,\tau)$ are expressed as a linear combination of power 
series of the form $q^{\lambda} \sum_{n=0}^\infty
\alpha_n (u)q^{n/\abs{g}}$.   
Moreover, both $S_i$ and $S^I$ satisfy the equations (\ref{m4.15}),
(\ref{m4.16}) and (\ref{m4.19}). 
Therefore, we shall investigate the power series with such properties
for a while. 

\begin{lem}
  Assume a formal power series $T(u,\tau ) = q^\lambda
  \sum_{n=0}^\infty \alpha_n(u) q^{n/\abs{g}}$ satisfies the
  conditions (\ref{m4.15}), (\ref{m4.16}) and (\ref{m4.19}). 
  Then the coefficient of the leading term $\alpha_0$ satisfies the
  following properties for $a\in V$ and $u\in U$ such that $ga=a$,
  $\phi(g)u=u$,  $ha=\xi a$, $\phi (h)u=\zeta u$ with 
  $\xi,\zeta\in\C :$ 
  \begin{align}
    &(i)\q \alpha_0(a*_gu)=\xi^{-1}\delta_{\xi\zeta,1}\alpha_0(u*_ga),
    \vsb\\
    &(ii)\q
       \alpha_0(u')=0 \q \text{for}\ u'\in O_g(U)+\dsum_{r\ne 0}U^r,
    \vsb\\ 
    & (iii)\q
      \alpha_0((\w*_g -c/24-\lambda )^N v)=0 \q \text{for}\
      \text{any}\ v\in U. 
  \end{align}
\end{lem}

\pf
First, we prove (ii).
Recall that $O_g(U)+\sum_{r\ne 0}U^r$ is generated by $w\in U^r$,
$r\ne 0$ and $b\circ_g v$ with $b\in V^r$, $v\in U^{\abs{g}-r}$, 
where $V^r=\{ x\in V \mid g x= e^{2\pii r/\abs{g}}x\}$ and 
$U^s=\{ y\in U \mid \phi (g) y=e^{2\pii s/\abs{g}}y\}$.
Since $w$ is contained in $O(g,h)$, we have $T(w,\tau)=0$ and hence 
$\alpha_0(w)=0$. 
If $r\ne 0$, then by Lemma \ref{qlem.9.2} (ii) we
get $\alpha_0(b\circ_g v)=0$.
Therefore we should show $\alpha_0(b\circ_g v)=0$ for $b\in V^0$ and
$v\in U^0$. If $hb=b$, then by Lemma \ref{qlem.9.2} (i) we have
$\alpha_0(b\circ_gv + \fr{1}{12}b_{[0]}v)=0$. 
On the other hand, by the definition of $O(g,h)$, we have 
$\alpha_0(b_{[0]}v)=0$ and hence $\alpha_0(b\circ_g v)=0$.
If $hb=\epsilon b$ with $\epsilon\ne 1$, then by Lemma \ref{qlem.9.2}
(ii) we obtain $\alpha_0(b\circ_g v)=0$.
Thus the assertion (ii) holds.

Let $a$, $u$, $\xi$ and $\zeta$ be as stated.
Then one can show the following equality.
\begin{equation}\label{above}
  a*_g u-u*_g a = a_{[0]}u.
\end{equation}
Suppose $\xi =\zeta =1$. In this case one has $\alpha_0(a_{[0]}u)=0$
and so that $\alpha_0(a*_g u)=\alpha_0(u*_g a)$.
If $\xi\ne 1$, then by Lemma \ref{qlem.9.2} (ii) the constant term
of $\sum_{k=0}^\infty Q_k(1,\xi,q)a_{[k-1]}u$ is equal to $-a*_g
u+\fr{1}{1-\xi}a_{[0]}u$ (note that in this case $a\not\in V^A$ so 
we don't know if $T(a_{[0]}u,\tau )=0$). 
Since $a\in V^0$ and $T$ vanishes on $\sum_{k=0}^\infty Q_k(1,\xi,\tau
) a_{[k-1]}u$ in this case, we can use the relation \eqref{above} 
to have 
$$
  \alpha_0(a*_g u)=\dfr{1}{1-\xi}(\alpha_0(a*_gu)-\alpha_0(u*_g a)).
$$
Solving this equation yields (i).

Consider (iii).
The constant term of $\tilde{\w}*_\tau v$ is equal to 
$\tilde{\w}*_g u - \fr{1}{2} \tilde{\w}_{[0]} v$ by Lemma
\ref{qlem.9.2}.
Since $T(\tilde{\w}_{[0]} v,\tau )=0$, the leading
term of 
$$
 \l(\tilde{\w}*_\tau -\fr{1}{\abs{g}} q_{\fr{1}{\abs{g}}}
  \fr{d}{dq_{\fr{1}{\abs{g}}}}\r) T(v,\tau )
$$ 
is equal to
$q^\lambda\alpha_0((\w*_g -c/24-\lambda )v)$. 
Since the operator $\ds \tilde{\w}*_\tau - \fr{1}{\abs{g}}
q_{\fr{1}{\abs{g}}} \fr{d}{dq_{\fr{1}{\abs{g}}}}$ will not decrease
the minimal degree of $T(v,\tau )$, the condition
$$
  \l(\tilde{\w}*_\tau -\fr{1}{\abs{g}} q_{\fr{1}{\abs{g}}}
  \fr{d}{dq_{\fr{1}{\abs{g}}}}\r)^N T(v,\tau ) =0
$$
implies
$\alpha_0((\w *_\tau -c/24-\lambda )^Nv)=0$.
Thus (iii) holds.
\qed
\vsv\\
By this lemma, the leading term $\alpha_0$ of the function $T$
with properties (\ref{m4.15}), (\ref{m4.16}) and (\ref{m4.19}) can be
seen as a linear function on $A_g(U)$ with the following $h$-twisted 
symmetric property: for $a\in A_g(V)$ and $u\in A_g(U)$ such that
$ha=\xi a$, $\phi (h)a=\zeta a$ with $\xi,\zeta\in\C$,
$$
  \alpha_0(a*_gu)=\xi^{-1}\delta_{\xi\zeta,1}\alpha_0(u*_ga) .
$$
We shall give an explicit description of such a linear function
$\alpha_0$.

\begin{prop}\label{mondai}
  Let $A$ be a finite-dimensional semisimple associative algebra over
  $\C$ and let $B$ be an $A$-bimodule.
  Let $h$ be an automorphism of an algebra $A$ of finite order.
  Assume that a representation $\phi : \lfm h\rfm \to \aut (B)$
  of a cyclic group $\lfm h\rfm$ on $B$ is given and it satisfies
  \begin{equation}\label{associative_condition}
    \phi (h) (a\cd b)=ha\cd \phi (h)b \q \text{and}\q 
    \phi (h) (b\cd a)=\phi (h)b\cd ha 
  \end{equation}
  for $a\in A$ and $b\in B$.
  Let $\w\in Z(A)$ such that $h\w =\w$. 
  Let $F$ be a linear function on $B$ with the following properties:
  \vsb\\
  (i)\ There exists a scalar $r\in \C$ and $N\gg 0$ such that
  \begin{equation}\label{center}
    F((\w -r)^N\cd b)=0 .
  \end{equation}
  (ii)\ For $ha=\alpha a$, $\phi (h)b=\beta b$ with
  $\alpha,\beta\in\C$,  
  \begin{equation}\label{twisted_function}
    F(a\cd b) =\alpha^{-1} \delta_{\alpha\beta,1}F(b\cd a) .
  \end{equation}
  Then $F$ can be expressed as
  \begin{equation}\label{top_trace}
    F(b)=\sum_j C_j\tr_{W^j}I_j(b)\psi_j (h)^{-1},
  \end{equation}
  where $(W^j,\pi_j)$, $\pi_j : A\to \End (W^j)$, is an $h$-stable left
  $A$-module on which 
  $\omega$ acts as a scalar $r$, $\psi_j (h)$ is a linear isomorphism
  on $W^j$ such that 
  $$
    \pi_j (ha)=\psi_j(h) \pi_j(a)\psi_j(h)^{-1},
  $$ 
  and $I_j$ is an element of $\hom_A(B\tensor_A W^j,W^j)$ such that
  $$
    \psi_j(h) \cd I_j(b)w= I_j\l(\phi (h)b\r)\cd \psi_j(h)w 
  $$
  for $a\in A$, $b\in B$ and $w\in W^j$.
\end{prop}

\pf
By the semisimplicity, $A$ decomposes into a direct sum of simple
two-sided ideals $A=\oplus_{i=1}^n A_i$.
As an $A$-bimodule, $B$ is also completely reducible.
Furthermore, the following lemma shows that $B$ decomposes into  
a direct sum of $\phi (h)$-invariant $A$-subbimodules.

\begin{lem}
  Under the assumption above, for every $\phi (h)$-invariant
  $A$-subbimodule $B'$ of $B$ there exists a $\phi (h)$-invariant
  $A$-subbimodule $B''$ such that  $B=B'\oplus B''$.
\end{lem}

\pf
Since $B$ is a semisimple $A$-bimodule, there exists an
$A$-subbimodule $B''$ such that $B=B'\oplus B''$. 
We show that we can choose a $\phi (h)$-invariant complement $B''$.
By the decomposition of $B$ above, every element $b\in B$ can be
written uniquely as $b=b_1+b_2$ for some $b_1\in B'$ and $b_2\in
B''$. 
Denote the projection map $b\mapsto b_2$ with respect to 
the decomposition above by $\sigma$. 
Clearly, $\sigma$ is an $A$-homomorphism. 
Using this projection we define a linear transformation $\rho$ on $B$
by 
$$
  \rho (b)=\fr{1}{\abs{h}}\sum_{i=0}^{\abs{h}-1}\phi (h)^{-i}
  \sigma \l( \phi (h)^ib\r) .
$$
Set $B'''=\rho B$. We claim that $B'''$ is a $\phi (h)$-invariant
$A$-subbimodule and $B=B'\oplus B'''$. 
For the $\phi (h)$-invariance, choose any $\rho (b)\in B'''$. Then we
have 
$$
\begin{array}{l}
  \phi (h)\rho (b)=\phi (h)\cd \dfr{1}{\abs{h}}
  \dsum_{i=0}^{\abs{h}-1}\phi (h)^{-i} \sigma\l(\phi (h)^i b\r)
  \vsb\\
  = \dfr{1}{\abs{h}} \dsum_{i=0}^{\abs{h}-1} \phi (h)^{-i+1}
    \sigma\l(\phi (h)^{i-1}\cd \phi (h)b\r)
  \vsb\\
  = \rho \l(\phi (h)b\r) .
\end{array}
$$
Thus $B'''$ is $\phi (h)$-invariant. Similarly we can show that
the left invariance $a\cd\rho (b)=\rho (a\cd b)$ and the right
invariance $\rho (b)\cd a=\rho (b\cd a)$ for all $a\in A$. 
Therefore, $B'''$ is an $A$-subbimodule.
Using $\rho^2=\rho$, one can show that $B=B'\oplus B'''$.
\qed
\vsb\\
By the preceding lemma, we may assume that $B$ doesn't have a proper
$\phi (h)$-invariant $A$-subbimodule. 
Since $B$ is a completely reducible $A$-bimodule,
we can take an irreducible $A$-subbimodule $B_1$ of $B$.
Clearly $\sum_{j=0}^{\abs{h}-1}\phi (h)^jB_1$ is a non-trivial $\phi
(h)$-invariant $A$-subbimodule of $B$, so one must have
$B=\sum_{j=0}^{\abs{h}-1}\phi (h)^jB_1$. 
Take a non-negative integer $t$ which is a maximal subject to the
condition that $\sum_{j=0}^{t-1}$ $\phi (h)^j$ $B_1$ $=$
$\oplus_{j=0}^{t-1}\phi (h)^jB_1$. 
Then by induction we see that $\phi (h)^lB_1\subset
\oplus_{j=0}^{t-1}\phi (h)^j B_1$ for all $l\geq 0$ and hence
$B=\oplus_{j=0}^{t-1}\phi (h)^jB_1$.
Since  $B_1$ is an irreducible $A$-bimodule, only one simple component
of $A$ gives a non-trivial left action on $B_1$.
Let us denote such a component by $A_1$. 
The assumption \eqref{twisted_function} is equivalent to the fact that 
\begin{equation}\label{twisted_function2}
  F\l( a\cd b\r) = F\l( b \cd h^{-1}a\r)\ \text{and}\ F(b)=0\ 
  \text{for}\ \phi (h)b=\lambda b,\ \lambda\ne 1.
\end{equation}
Moreover, one can equivalently replace the condition above by
\begin{equation}\label{twisted_function3}
  F\l(a\cd b\r) =F\l(\phi (h)^{-1}b\cd h^{-2}a\r) .
\end{equation}
Let $1_{A_1}$ be the unit element of the simple ring $A_1$. For every
$b_1\in B_1$, the condition \eqref{twisted_function2} implies that
$$
  F(b_1)= F\l( 1_{A_1}\cd b_1\r) =F\l( b_1\cd h^{-1}1_{A_1}\r) 
  = F\l( b_1\cd 1_{h^{-1}A_1}\r) .
$$
So if $B_1$ is not a non-trivial right $h^{-1}A_1$-module, then
$F(B_1)=0$, and by \eqref{twisted_function3} we arrive at $F(B)=0$.
Therefore, we may assume that $B_1$ is a $(A_1, h^{-1}A_1)$-bimodule.
Since $h$ permutes the simple components of $A$, each irreducible
component $\phi (h)^jB_1$ of $B$ is an irreducible
$(h^jA_1,h^{j-1}A_1)$-bimodule. 
So we may assume that $h$ acts transitively on $A$ as a permutation of
its simple components. 
This means, we may assume that $A=A_1\oplus hA_1\oplus \cds \oplus 
h^{s-1}A_1$ and $h$ maps $h^{s-1}A_1$ onto $A_1$.
In particular, $h^s$ induces an automorphism on each simple component
$h^iA_1$ of $A$ for $0\leq i<s$. 
Since $\phi (h)^jB_1$ is an irreducible $(h^j A_1,h^{j-1}
A_1)$-bimodule, we note that $s\leq t$.
Furthermore, $\phi (h)B=B$ forces $\phi (h)^tB_1$ to be an
irreducible $(A_1,h^{-1}A_1)$-bimodule so that we obtain $s|t$.
As $B=\oplus_{j=0}^{t-1}\phi (h)^jB_1$, one can find that the
number of isotypical components of $B$ as an $A$-bimodule is $s$.
Therefore, $B$ decomposes into a direct sum of its isotypical
components as  
$$
  \displaystyle
  B=\bigoplus_{k=0}^{s-1}B^{(k)},\q \text{where}\ 
  B^{(k)}=\bigoplus_{i=0}^{t/s-1}\phi (h)^{si+k} B_1.
$$
Hence, $\phi (h)$ induces a cyclic permutation on the set of
isotypical components as $\phi (h) :B^{(k)}\to B^{(k+1)}$
and thus $\phi (h)^s$ induces a linear transformation on each
isotypical component of $B$.  

By definition, $B^{(0)}$ is an $(A_1,h^{-1}A_1)$-bimodule.
We twist the right action of $A$ on $B^{(0)}$ by $h$ to obtain an
$A_1$-bimodule structure on $B^{(0)}$. Define a new right action of
$a\in A$ by
$$
  b*a:=b\cd (h^{-1}a)\q \text{for}\  b\in B^{(0)}.
$$
Under the action above, $B^{(0)}$ becomes an $A_1$-bimodule and
$B^{(0)}$ decomposes into a direct sum of $t/s$ irreducible
$A_1$-subbimodules as $B^{(0)}=\oplus_{j=0}^{t/s-1} B_{0j}$. 
Since each irreducible component $B_{0j}$ is isomorphic to $A_1$ 
as an $A_1$-bimodule, there exist $A_1$-isomorphisms $f_j: A_1
\simto B_{0j}$ for $0\leq j\leq t/s-1$. 
Using these $f_j$'s, we define linear functions $G_j:A_1\to \C$ by
$$
  G_j: a\in A_1 \mapsto F\l( f_j(a_1)\r) \in \C ,
  \q 0\leq j\leq t/s-1.
$$
Then for $a_1,a_2\in A_1$, we have 
$$
\begin{array}{l}
  G_j(a_1a_2)=F\l( f_j(a_1a_2)\r) = F\l( a_1\cd f_j(a_2)\r)
  = F\l( f_j(a_2)\cd h^{-1}a_1\r) 
  \vsb\\
  = F\l( f_j(a_2)*a_1\r) = F\l( f_j(a_2\cd a_1)\r)
  = G_j(a_2a_1).
\end{array}
$$
Therefore, the linear functions $G_j$ are trace functions. 
Let $W$ be an irreducible left $A_1$-module. Then $G_j$ can be 
written as $G_j(a)=C_j\, \tr_W(a)$ with some $C_j\in\C$
and the linear function $F$ restricted on $B_{0j}$ can be described as 
\begin{equation}\label{trace}
  F(b_j)=C_j\, \tr_W f_j^{-1}(b_j)\q (C_j\in\C ) 
\end{equation}
for $b\in B_{0j}$.
To get a desired description of $F$, we have to extend $W$ to be an
$A$-module. 
Based on a left $A_1$-module structure on $W$, we construct left
$h^kA_1$-modules $h^k\circ W$, $0\leq k\leq s-1$. 
As a vector space, we set $h^k\circ W=W$ and we denote the elements of
$h^k\circ W$ formally by $h^kw$, where $w\in W$.  
We introduce an action of $h^kA_1$ on $h^k\circ W$ by
$h^ka\cd h^kw:=h^k(a\cd w)$. By this definition,  $h^k\circ W$ becomes 
an irreducible left $h^kA_1$-module. Set $\bar{W} = \oplus_{k=0}^{s-1}
h^k \circ W$. Then $\bar{W}$ naturally carries a left $A$-module
structure.  
Recall that $A=\oplus_{k=0}^{s-1}h^kA_1$ and $h^s$ is an automorphism
of each simple component of $A$.  
By Skolem-Noether theorem, an automorphism $h^s$ of a simple ring 
$A_1$ is realizable by an inner automorphism of $A_1$.
That is, there exists an element $\gamma\in A_1^*$ 
such that $h^sa_1=\gamma a_1\gamma^{-1}$ for all $a_1\in A_1$. 
Then the actions of $h^s$ on $h^kA_1$, $0\leq
k\leq s-1$, are given by inner automorphisms defined by
$h^k\gamma\in (h^kA_1)^*$.
Furthermore, an element $\gamma$ induces a linear isomorphism $\psi$
of $\bar{W}$ by the following way.
For $0\leq k\leq s-2$, we set $\psi ( h^kw) =h^{k+1}w$ and
we define $\psi (h^{s-1}w)= h^0\gamma w$ for $k=s-1$. 
It is easy to check that  $\psi$ is a linear isomorphism on $\bar{W}$. 
This automorphism is characterized by the following property:
$$
  \psi (a\cd w)=ha\cd \psi (w),\ a\in A,\ w\in \bar{W}.
$$ 
This means that $\bar{W}$ is an $h$-stable $A$-module with an
$h$-stabilizing automorphism $\bar{\psi}$.

Next we shall define a linear map $I_0:B^{(0)}\times \bar{W}\to
\bar{W}$. For $b_j\in B_{0j}$ and $w\in \bar{W}$, we set
$$
  I_0: (b_j, w)  \mapsto C_jf_j^{-1}(b_j)\cd \psi (w)
$$
and extend bilinearly on $B^{(0)}\times \bar{W}$. 
By this definition, we can show that $I_0$ defines an element of
$\hom_A(B^{(0)}\tensor_A \bar{W}, \bar{W})$.
By \eqref{trace}, we note that 
$$
  F(b)= \tr_{\bar{W}}I_0(b)\psi^{-1}
$$
for all $b\in B^{(0)}$.
Up to now, we have determined an explicit description of $F$ 
restricted on $B^{(0)}$. 
By \eqref{twisted_function3}, a description of $F$ on
$B^{(0)}$ completely determines that on $B$.
Set
$$
  \tilde{I}_0(b):=\dfr{s}{\abs{h}} \sum_{i=0}^{\abs{h}/s-1} 
    \psi^{-is} I_0\l( \phi (h)^{is}b\r)\psi^{is} .
$$
Since $\phi (h)^s$ is an automorphism on $B^{(0)}$ and $I_0\in \hom_A
(B \tensor_A \bar{W}, \bar{W})$, using \eqref{associative_condition}
we can show that $I_0'(b)=\psi^{-s}I_0\l(\phi (h)^s b \r)\psi^s$ is
also an element of $\hom_A$ $(B^{(0)}\tensor_A \bar{W},\bar{W})$. 
Furthermore, it follows from \eqref{twisted_function3}
that $F\l(\phi (h)^sb\r) =F\l( b\r)$, and thus $\tr_{\bar{W}} I_0'(b) 
\psi^{-1} = \tr_{\bar{W}} I_0(b)\psi^{-1} =F(b)$ for all $b\in
B^{(0)}$.  
Therefore, we have 
$$
  \tr_{\bar{W}} \tilde{I}_0(b) \psi^{-1} 
  = \tr_{\bar{W}}I_0(b)\psi^{-1}=F(b)
$$ 
for all $b\in B^{(0)}$. 
It is clear that  $\psi^{-s}\tilde{I_0}\l(\phi (h)^s b\r)\psi^s
=\tilde{I_0}(b)$.  
Hence, by replacing  $I_0$ by $\tilde{I}_0$, we may assume 
\begin{equation}\label{normalize}
  \psi^{-s}I_0\l( \phi(h)^s b \r) \psi^s = I_0(b).
\end{equation}

As we did before, by twisting the right action by $h^k$, we can
introduce an $h^kA_1$-bimodule structure on each isotypical component 
$B^{(k)}$, $0\leq k\leq t/s-1$. 
Since $B^{(k)}=\phi (h)^k B^{(0)}$, $B^{(k)}$ decomposes
into a direct sum of irreducible $h^kA_1$-bimodules as
\begin{equation}\label{decomp} 
  B^{(k)}=\bigoplus_{j=0}^{t/s-1} \phi (h)^k B_{0j}.
\end{equation}
On this decomposition, we define $I_k:B^{(k)}\times \bar{W}\to
\bar{W}$ as follows.
For $\phi (h)^kb_j$, $0\leq k\leq s-1$, $b_j\in
B_{0j}$ and $w\in\bar{W}$, we set
$$
  I_k\l(\phi (h)^kb_j\r) w:=\psi^kI_0\l(b_j\r)\psi^{-k}w
$$
and extend linearly on $B^{(k)}$. It follows from the definition that 
$I_k$ $\in$ $\hom_A$ $(B^{(k)}\tensor_A\bar{W}$ $,\bar{W})$. 
Furthermore, by \eqref{twisted_function3}, we see that $F(b) =
\tr_{\bar{W}} I_k(b) \psi^{-1}$ for all $b\in B^{(k)}$.
By extending the defining region of $I_k$ on the whole of $B$ naturally
and setting $\bar{I}:=I_0+I_1+\cds +I_{s-1}$, we see that
$\bar{I}\in \hom_A(B \tensor_A \bar{W},\bar{W})$ and
$F(b)=\tr_{\bar{W}}\bar{I}(b)\psi^{-1}$ for all $b\in B$.
Moreover, by \eqref{normalize}, $\psi\cd \bar{I}(b)w =
\bar{I}\l(\phi (h)b\r)\cd \psi (w)$ for all $b\in B$ and $w\in
\bar{W}$. Therefore $\bar{I}$ and $\psi$ satisfy all conditions we
required. Finally, since $\w\in Z(A)$, it is clear that $\w$ acts on
$\bar{W}$ as a scalar, which should be $r$ by the assumption. 
\qed

Now we are ready to apply the preceding proposition to $\alpha_0$ by
setting $A=A_g(V)$, $B=A_g(U)$ and $F=\alpha_0$. 
By the previous proposition and Theorem \ref{extended Li}, we can
construct intertwining operators $I_k(\cd\, ,z)$ among fundamental
$h$-stable $g$-twisted $V$-modules such that the traces of the zero
mode on the top modules are given by $\alpha_0$. 
Furthermore, the difference between the actions of $h$-stabilizing
automorphisms on the top levels of $h$-stable $V$-modules and those of
$h$-stabilizing automorphisms $\psi$ given in the previous proposition
are only scalar multiples. 
Hence, an $h$-stabilizing automorphism $\psi$ on a top level gives
rise to an $h$-stabilizing automorphism on $V$-module.
Thus we can obtain fundamental $h$-stable $\phi(g)$-twisted
intertwining operators $I_k(\cd\, ,z)$ such that $\alpha_0$ is a linear
combination of the leading terms of $S^{I_k}(\cd,\tau)$.
Therefore, we have the following statement.  

\begin{prop}
  Let $\alpha_0$ be the coefficient of leading term of the
  power series $T(u,\tau )$ $=$ $q^\lambda$ $\sum_{n=0}^\infty$
  $\alpha_n (u)$ $q^{n/\abs{g}}$ which satisfies \eqref{m4.15},
  \eqref{m4.16} and \eqref{m4.19}. 
  Then there exist fundamental $h$-stable $g$-twisted $V$-modules 
  $(W_g^k,\psi_k (h))$ with minimum $L(0)$-weights $\lambda +c/24$ and 
  fundamental $h$-stable $\phi (g)$-twisted 
  intertwining operators $I_k(\cd,z)$ of type $U\times W_g^k\to W_g^k$
  such that  
  $$
    \alpha_0 (u)=\sum_k C_k\ \tr_{|_{W_g^k}} o^{I_k}(u)\psi_k(h)^{-1}.
  $$
\end{prop}

Thus we have

\begin{prop}\label{mlem.4.13}
   Suppose $U$ is $C_{[2,0]}^A$-cofinite.
   Assume a formal power series $T(u,\tau )=q^\lambda\sum_{n=0}^\infty
   \alpha_{n/\abs{g}}(u) q^{n/\abs{g}}$ is given and satisfies
   (\ref{m4.15}),  (\ref{m4.16}) and (\ref{m4.19}).
   Then there exist fundamental $h$-stable $g$-twisted $V$-modules 
   $(W_g^k,\psi_k(h))$ and $h$-stable $\phi (g)$-twisted intertwining 
   operators $I_k(\cd\, ,z)$ of type $U\times W_g^k\to W_g^k$ such
   that $T(u,\tau )$ can be expressed as a linear combination of the
   trace functions $S^{I_k}(u,\tau )$. 
\end{prop} 

\pf
By Proposition \ref{mondai}, there exist $h$-stable $A_g(V)$-modules
$(W_g^k (0)$, $\psi_k(h))$, $1\leq k\leq m$, on which 
$\omega$ acts as a scalar $\lambda +c/24$ and linear maps
$$
  I_k \in
  \hom_{A_g(V)} ( A_g(U) \tensor_{A_g(V)} W_g^k(0) \circ h,W_g^k(0)
  ),\q 1\leq k\leq m
$$ 
such that
$$
  \alpha_0(u)=\sum_{j=1}^mC_j\, \tr_{|_{W_g^j(0)}}I_j(u)\psi_j(h)^{-1}.
$$
By Theorem \ref{extended Li}, there exist fundamental $h$-stable
$g$-twisted $V$-modules $(W_g^k, \psi_k(h))$ and fundamental
$h$-stable $\phi (g)$-twisted intertwining operators 
$\bar{I_k}(\cd\,,z)$ of type $U\times W_g^k\to W_g^k$ such that
all the minimum $L(0)$-weights of $W_g^k$ are equal to $\lambda +c/24$ 
and whose top levels are isomorphic to $W_g^j(0)$ as $A_g(V)$-modules
for $1\leq k\leq m$.   
Then the function $T(u,\tau )-\sum_j S^{I_j}(u,\tau )$ also
satisfies the conditions (\ref{m4.15}), (\ref{m4.16}) and
(\ref{m4.19}) and the degree of its leading term is strictly greater
than $\lambda$. 
If $T(u,\tau )-\sum_j S^{I_j}(u,\tau )$ is not equal to $0$, then we
can repeat the same argument on it. 
But by the rationality of $V$, there exist only finitely
many inequivalent irreducible $V$-modules. 
So our procedure must stop in finite steps. 
Hence, $T(u)$ can be expressed as a linear combination of
$S^{I_j}(u)$'s as desired.  
\qed

\begin{prop}
  $S_i(\cd,\tau )=0$ for $i>0$.
\end{prop}

\pf
Assume $p\geq 1$. By Proposition \ref{mlem.4.13}, $S_{p-1}(\cd,\tau )$ 
can be written as a linear combination of the trace functions.
Since the trace functions satisfy the condition 
$$
  S^I( \tilde{\w}*_\tau u,\tau ) 
  =\fr{1}{2\pii} \fr{d}{d\tau}S^I(u,\tau ),
$$
we have 
$$
  S_{p-1}(\tilde{\w}*_{\tau}u,\tau ) 
  =\fr{1}{2\pii} \fr{d}{d\tau} S_{p-1} (u,\tau ).
$$
On the other hand, by \eqref{m4.18}, we have
$$
  S_{p-1}(\tilde{\w}*_\tau u, \tau ) 
  = p S_p(u,\tau )+ \fr{1}{2\pii} \fr{d}{d\tau} S_{p-1}(u,\tau ).
$$ 
Therefore we get $pS_p(u,\tau )=0$, a contradiction. 
\qed
\vsb\\
Summarizing up to everything, we arrive at the following main
theorem. 
This is a new extension of Zhu theory which generalize both the
orbifold case in \cite{DLM2} and the intertwining operators case in
\cite{M3}.  

\begin{thm}\label{conclusion}
  Let $V$ be a $g$-rational VOA and $U$ an $C_{[2,0]}^A$-cofinite 
  $A$-stable untwisted $V$-module with a stabilizing automorphism
  $\phi$.
  Let $\{ (W_g^1,\psi_1(h)),\dots,(W_g^m,\psi_m(h))\}$ be the set of
  all inequivalent fundamental $h$-stable $g$-twisted $V$-modules 
  and let $\{ I_{kj}(\cd,z) \mid j=1,\dots ,r_k\}$ be a basis
  of fundamental $h$-stable $\phi (g)$-twisted intertwining operators
  of type $U\times W_g^k\to W_g^k$ for each $1\leq k\leq m$. 
  Then the space of 1-point functions $\mathcal{C}_1(U;\phi;(g,h))$ is
  spanned by the trace functions 
  $$
    S^{I_{kj}}(u,\tau )=\tr_{W_g^k} z^{\wt u}I_{kj}(u,z)
    \psi_k(h)^{-1} q^{L(0)-c/24}
  $$
  for $1\leq j\leq r_k$, $1\leq k\leq m$.
  In particular, the linear space spanned by the trace functions
  $$
    \l\la\,  S^{I_{kj}}(\cd,\tau)\ \Big|\ 1\leq j\leq r_k, 1\leq
    k\leq m \, \r\ra
  $$ 
  is invariant under the modular group $\Gamma (g,h)$,
  where its action is defined by \eqref {q5.10}.   
\end{thm}

\begin{center}
  \gs{Acknowledgments}
\end{center}
\begin{quote}
  The author wishes to thank Professor Masahiko Miyamoto for his
  insightful advice. 
  He also thanks to Professor Ching Hung Lam for reading this
  manuscript and his helpful comments. 
\end{quote}

\small
\setlength{\baselineskip}{12pt}


\begin{thebibliography}{hoge99}

\bibitem[Ab]{Ab}
  T. Abe, Fusion rules for the free bosonic orbifold vertex operator
  algebra, J. Algebra 229, 333-374 (2000)

\bibitem[AM]{AM}
  G. Anderson and G. Moore, Rationality in conformal field theory,
  Commu. Math. Phys. 117 (1988), 441-450.

\bibitem[Bo]{Bo}
  R.E. Borcherds, Vertex operator algebras, Kac-Moody algebras and the 
  Monster, Proc. Natl. Acad. Sci. USA. 83 (1986), 3026

\bibitem[DLM1]{DLM1}
  C. Dong, H. Li, G. Mason, Twisted representations of vertex operator
  algebras, Math. Ann. 310 (1998) 571-600. 

\bibitem[DLM2]{DLM2}
  C. Dong, H. Li, G. Mason, Modular-invariance of trace functions in
  orbifold theory, Commu. Math. Phys. 214, 1-56.

\bibitem[DLM3]{DLM3}
  C. Dong, H. Li, G. Mason, Simple currents and extensions of vertex
  operator algebras, Commu. Math. Phys. 180, 671-707.

\bibitem[DM]{DM}
  C. Dong, G. Mason, On Quantum Galois Theory, Duke Math. J. vol.86,
  No.2 (1997) 305-321.

\bibitem[FHL]{FHL}
  I. Frenkel, Y.-Z. Huang and J. Lepowsky, On axiomatic approaches to
  vertex operator algebras and modules, Memoirs Amer. Math. Soc. 104,
  1993

\bibitem[FLM]{FLM}
  I.B. Frenkel, J. Lepowsky, A. Meurman, Vertex Operator Algebras and
  the Monster, Academic Press, New York, 1988.

\bibitem[FZ]{FZ}
  I.B. Frenkel, Y. Zhu, Vertex operator algebras associated to
  representation of affine and Virasoro algebras, Duke Mathematical
  Journal 66 (1992), 123-168.

\bibitem[La]{La}
  S. Lang, Elliptic Functions, Springer-Verlag, 1987.

\bibitem[Li1]{Li1}
  H. Li,  Local systems of twisted vertex operators, vertex operator
  superalgebras and twisted modules, Cont. Math. 193 (1996) 203-236 

\bibitem[Li2]{Li2}
  H. Li, Determining fusion rules by $A(V)$-modules and bimodules, 
  J. Algebra 212, 515-556.

\bibitem[Li3]{Li3}
  H. Li, An analogue of the Hom functor and a generalized nuclear
  democracy theorem, Duke Math. J. 93 (1998), 73-114. 

\bibitem[M1]{M1}
  M. Miyamoto, A modular invariance on the theta functions defined on
  vertex operator algebras, Duke Math. J. Vol. 101 (2000), 221-236.

\bibitem[M2]{M2}
  M. Miyamoto, Modular invariance of trace functions on VOA in many
  variables, preprint.

\bibitem[M3]{M3}
  M. Miyamoto, Intertwining operators and modular invariance, 
   q-alg/0010180.

\bibitem[X]{X}
  X. Xu, Intertwining operators for twisted modules of a colored
  vertex operator superalgebra, J. Algebra 175, 241-273.

\bibitem[Z]{Z}
  Y.Zhu, Modular invariance of characters of vertex operator algebras, 
  J. Amer. Math. Soc. 9 (1996), 237-302.

\end{thebibliography}
\end{document}